\documentstyle{amsppt}
\magnification1200
\pagewidth{6.5 true in}
\pageheight{9.25 true in}
\NoBlackBoxes

\def\phi{\varphi}

\def\lam{\lambda}
\topmatter
\title Omega results for the divisor and circle problems
\endtitle
\author K. Soundararajan
\endauthor
\address{Department of Mathematics, University of Michigan, Ann Arbor,
Michigan 48109, USA} \endaddress
\email{ksound{\@}umich.edu} \endemail
\thanks{The author is partially supported by the
National Science Foundation and by the American
Institute of Mathematics (AIM).}
\endthanks
\endtopmatter

\document

\head 1. Introduction \endhead

\noindent Let $d(n)$ denote the number of divisors of $n$ and $r(n)$ the
number of ways of writing $n$ as the sum of two integer squares.  Let
$\Delta(x)$ and $P(x)$ denote the remainder terms in the aymptotic
formulae $\sum_{n\le x} d(n) = x\log x + (2\gamma -1) x +\Delta(x)$
and $\sum_{n\le x} r(n) = \pi x + P(x)$.  In 1916 G.H. Hardy [4]
showed that
$$
\Delta(x) =
\cases
\Omega_+((x\log x)^{\frac 14} \log_2 x)\\
\Omega_-(x^{\frac 14})\\
\endcases
$$
and that
$$
P(x) = \cases
\Omega_- ((x\log x)^{\frac 14}) \\
\Omega_+ (x^{\frac 14}).\\
\endcases
$$
Here and throughout $\log_j$ denotes the $j$-th iterated logarithm, so
that $\log_2 =\log \log$, $\log_3 =\log \log \log $ and so on.  Recall
that for a real valued function $f$ and a positive function $g$ the symbol
$f=\Omega(g)$ means that $\lim \sup_{x\to \infty} |f(x)|/g(x) >0$.
We write $f=\Omega_+(g)$ if $\lim \sup_{x\to \infty} f(x)/g(x) >0$, and
$f=\Omega_-(g)$ if $\lim \inf_{x\to \infty} f(x)/g(x) <0$.  Lastly
$f=\Omega_\pm (g)$ means that $f=\Omega_+(g)$ and also $f=\Omega_-(g)$.

Since
Hardy, gradual progress had been made on the $\Omega_-$ result
for $\Delta$ and the $\Omega_+$ result for $P$ culminating
in the work of K. Corr{\' a}di and I. K{\' a}tai [1] who showed that
for a positive constant $c$
$$
\Delta(x) = \Omega_-\Big( x^{\frac 14}
\exp(c (\log_2 x)^{\frac 14} (\log_3 x)^{-\frac 34})\Big)
$$
and a similar $\Omega_+$ result for $P(x)$.  In 1981 J.L. Hafner [2]
obtained the first improvements on the $\Omega_+$ result for $\Delta$
and the $\Omega_-$ result for $P$.   He showed that for some
positive constants $A$ and $B$,
$\Delta(x)=\Omega_+ ( (x\log x)^{\frac 14} (\log_2 x)^{(3+2\log 2)/4}
\exp(-A \sqrt{\log_3 x}))$ and
$P(x) = \Omega_- ((x\log x)^{\frac 14} (\log_2 x)^{(\log 2)/4}
\exp(-B\sqrt{\log_3 x}))$.
Hafner observed that these results represented the limit of
his method and that A. Selberg (unpublished) had obtained similar
bounds.   In this note we refine Hafner's results and show
that the magnitudes of $\Delta(x)$ and $P(x)$ can be larger
than the values given above.  However, unlike Hafner's result,
we cannot determine the sign of the large values we exhibit.

\proclaim{Theorem 1}  We have
$$
\Delta(x) = \Omega\Big( (x\log x)^{\frac 14} (\log_2 x)^{\frac{3}{4}
(2^{4/3}-1)} (\log_3 x)^{-\frac 58}\Big),
$$
and
$$
P(x)= \Omega\Big((x\log x)^{\frac 14} (\log_2 x)^{\frac  34(2^{1/3}-1)}
(\log_3 x)^{-\frac 58}\Big).
$$
\endproclaim

Note that $\frac 34 (2^{4/3}-1)= 1.1398\ldots$ while
$(3+2\log 2)/4= 1.0965\ldots$; also $\frac  34(2^{1/3}-1) = 0.1949\ldots$
while $(\log 2)/4 = 0.1732\ldots$.

Our method also
applies to the remainder term in the $k$-divisor problem (also called
the Piltz divisor problem).  Let $k\ge 2$ be an integer
and let $d_k(n)$ denote the number
of ways of expressing $n$ as a product of $k$ factors.  Let
$\Delta_k(x)$ denote the remainder term in the asymptotic
formula for $\sum_{n\le x} d_k(n)$;  that is,
$$
\sum_{n\le x} d_k(n) = \mathop{\text{Res}}_{s=1}\  \zeta(s)^k \frac{x^s}{s}
+ \Delta_k(x).
$$
G. Szeg{\" o} and A. Walfisz [7, 8] showed that $\Delta_k(x)
= \Omega^*( (x\log x)^{(k-1)/(2k)} (\log_2 x)^{k-1})$
where $\Omega^*= \Omega_+$ if $k=2$, $3$ and $\Omega^* =\Omega_\pm$
if $k\ge 4$.  Hafner [3] improved this to
$$
\Delta_k(x)= \Omega^*\Big( (x\log x)^{\frac{k-1}{2k}}
(\log_2 x)^{\frac{(k-1)}{2k}(k\log k-k+1)+k-1} \exp(-A_k\sqrt{\log_3 x})\Big)
$$
for some positive constant $A_k$.
We exhibit larger values of $|\Delta_k(x)|$ but as in
Theorem 1 we cannot control the sign of these values (except
when $k\equiv 3\pmod 4$).

\proclaim{Theorem 2}  With notations as above
$$
\Delta_k(x)= \Omega\Big( (x\log x)^{\frac{k-1}{2k}}
(\log_2 x)^{\frac{k+1}{2k} (k^{2k/(k+1)} -1)} (\log_3 x)^{-\frac 12
-\frac{k-1}{4k}}
\Big).
$$
The above estimate holds with $\Omega_+$ in place of $\Omega$ if
$k\equiv 3\pmod 8$, and with $\Omega_-$ in place of $\Omega$ if
$k\equiv 7\pmod 8$.
\endproclaim

For large $k$ the exponent of $\log_2 x$ in our
result is $\sim k^2/2$ while that in Hafner's is $\sim (k\log k)/2$.

We now describe our method, using
$\Delta(x)$ for illustration.  One knows that $\Delta(x^2)$ is
given by the conditionally convergent series
$$
\frac{x^{\frac 12}}{\pi\sqrt{2}}
\sum_{n=1}^{\infty} \frac{d(n)}{n^{\frac 34}} \cos (4\pi\sqrt{n}x -\pi/4).
$$
By smoothing a little, one may restrict the sum above to the terms $n\le N$
weighted appropriately, and it suffices (roughly speaking)
to give omega results for the truncated series
$\sum_{n\le N} d(n)n^{-\frac 34}\cos (4\pi\sqrt{n}x -\pi/4)$.
Let ${\Cal M}$ denote a set of $M$ positive integers.
By Dirichlet's Theorem on diophantine approximation we may find
$x\in [X,6^M X]$ such that
$\parallel 2\sqrt{m}x\parallel \le 1/6$ for each $m\in {\Cal M}$.
\footnote{Here $\parallel \cdot\parallel$ denotes the distance
from the nearest integer.}
If we select ${\Cal M}$ to be the first $M$ integers, and
take $M= [\log X]=N$ then we obtain Hardy's omega result.  Hafner exploits
the uneven distribution of $d(n)$ by selecting ${\Cal M}$ such that
$d(m)$ is large for $m\in {\Cal M}$.  To ensure
that the terms $n\le N$, $n\notin {\Cal M}$ do not cancel the
contribution of the terms $m\le N$, $m\in {\Cal M}$, Hafner imposes the
restriction $\sum\Sb n\le N, n\notin {\Cal M}\endSb d(n)n^{-\frac 34}
=o(N^{\frac 14}\log N)$.  Optimizing this argument leads to his $\Omega_+$
result.  We argue instead as follows: For an integer parameter $L$, we
first find $x\in [X,(6L)^M X]$ such that $\parallel 2\sqrt m x\parallel
\le 1/(6L)$ for each $m\in {\Cal M}$.  Then for each of the $L$ points
$\ell x$ ($1\le \ell \le L$) we see that the terms $m\le M$, $m\in {\Cal M}$
pull in the same direction.  We then show that for one of these points
the contribution of the terms $n\le N$, $n\notin {\Cal M}$ is
not too destructive.  The effect is essentially to eliminate
Hafner's restriction, and this accounts for our improvement.
Our argument really works for
$\sum_{n\le N} d(n)n^{-\frac 34} \cos(4\pi\sqrt{n}x)$, so that
it is first necessary to remove the phase $-\pi/4$.  It is
in this step that we lose knowledge of the sign of the
large values we exhibit.

From our remarks above the ideal omega result for $\Delta(x)$
seems the following.  Arrange the sequence $d(n)n^{-\frac 34}$ in
descending order, and let $S(M)$ denote the sum of the first
$M$ largest values.  Then $\Delta(x) = \Omega(x^{\frac 14} S(\log x))$.
One can show that
$S(M) = M^{\frac 14} (\log M)^{\frac 34(2^{4/3}-1) +o(1)}$;
thus Theorem 1 essentially obtains this ideal omega result.

We may model \footnote{This is not entirely accurate since the terms
at $n$ and $nm^2$ are obviously correlated.  With this caveat
the model is plausible, see D.R. Heath-Brown [5].}
$\pi\sqrt{2}\Delta(x)x^{-\frac 14}$ by a
random trigonometric series
$\sum_{n=1}^{\infty} d(n)n^{-\frac 34} \cos(X_n)$ where
the $X_n$ are independent random variables uniformly distributed on
$[0,2\pi)$.  The work of H.L. Montgomery and A.M. Odlyzko [6]
provides estimates for the probability of large values attained
by this trigonometric series.  This suggests that the omega result
obtained in Theorem 1 represents the true maximal order of $\Delta(x)$ up to
$(\log_2 x)^{o(1)}$.

\head 2. The key Lemma \endhead

\noindent Let $f(1)$, $f(2)$, $\ldots$ be a sequence of non-negative real
numbers and $0 \le \lam_1 \le \lam_2 \le \ldots$ be a non-decreasing
sequence of non-negative real numbers.  We suppose that
$\sum_{n=1}^{\infty} f(n) < \infty$ and consider the trigonometric
series
$$
F(x) := \sum_{n=1}^{\infty} f(n) \cos(2\pi \lam_n x +\beta)
$$
where $\beta \in {\Bbb R}$.

\proclaim{Lemma 3} Let $L\ge 2$ and $N\ge 1$ be integers.
Let ${\Cal M}$ be a set of integers such that
$\lam_m \in [\frac{\lam_{N}}2, \frac{3\lam_N}{2}]$ for each $m\in {\Cal M}$.
For any $X\ge 2$ there exists a point
$x\in [X/2,(6L)^{M+1}X]$ such
that
$$
|F(x)|  \ge  \frac 18 \sum_{m\in {\Cal M}} f(m)  -
\frac{1}{L-1}\sum\Sb n\\ \lam_n \le 2\lam_N\endSb f(n) -
\frac{4}{\pi^2 X\lam_N} \sum_{n} f(n).
\tag{1}
$$
If $\beta \equiv 0 \pmod{2\pi}$ then
there is a point $x\in [X/2,(6L)^{M+1}X]$
such that
$$
F(x) \ge \frac 18 \sum_{m\in {\Cal M}} f(m)  -
\frac{1}{L-1}\sum\Sb n\\ \lam_n \le 2\lam_N\endSb f(n) -
\frac{2}{\pi^2 X\lam_N} \sum_{n} f(n).
\tag{2}
$$
If $\beta \equiv \pi \pmod{2\pi}$ then the conclusion (2)
holds with $-F(x)$ in place of $F(x)$.
\endproclaim
\demo{Proof} Let
$K(u) = (\frac{\sin(\pi u)}{\pi u})^2$ be Fejer's
kernel and recall that $\int_{-\infty}^{\infty} K(u) e(-uy)
du = \max(0, 1-|y|)=: k(y)$,
say.  Consider
$$
\align
\int_{-\infty}^{\infty} &  \lam_N K(\lam_N u) e(-\lam_N u)F(x+u) du\\
&=\frac 12 \sum_{n} f(n) \int_{-\infty}^{\infty} \lam_N
K(\lam_Nu) e(-\lam_N u)
\Big( e^{i\beta} e(\lam_n (x+u)) + e^{-i\beta} e(-\lam_n (x+u))\Big) du \\
&= \frac {e^{i\beta}}2 \sum_{n} f(n) e(\lam_n x)
k\Big(\frac{\lam_N-\lam_n}{\lam_N}\Big),\\
\endalign
$$
since $k((\lam_N+\lam_n)/\lam_N)=0$.  Setting
$$
F_1(x) = \frac 12\sum_{n} f(n) \cos(2\pi \lam_n x)
k\Big(\frac{\lam_N-\lam_n}{\lam_N}\Big),
$$
we deduce that
$$
\align
F_1(x) &\le
\int_{-\infty}^{\infty} \lam_N K(\lam_N u) |F(x+u)| du \\
&\le \int_{-X/2}^{X/2} \lam_N K(\lam_N u) |F(x+u)| du
+ \int_{|u|>X/2} \frac{1}{\pi^2 \lam_N u^2} \sum_{n} f(n) du
\\
&\le \max_{u\in [-X/2,X/2]} |F(x+u)| +\frac{4}{\pi^2 X \lam_N}\sum_{n} f(n).
\tag{3}
\\
\endalign
$$

By Dirichlet's Theorem (see for example \S 8.2 of [10]), for any $X\ge 2$
there exists a point $x_0$ in $[X,(6L)^MX]$ such that
$\parallel \lam_m x_0 \parallel \le 1/(6L)$ for each $m \in {\Cal M}$.
Consider
$$
\sum_{\ell = -L}^{L} k\Big( \frac{\ell}{L}\Big) F_1(\ell x_0)
= \frac{1}{2} \sum_{n} f(n) k\Big(\frac{\lam_N-\lam_n}{\lam_N}\Big)
\sum_{\ell = -L}^{L} k\Big(\frac {\ell}{L}\Big) \cos (2\pi \lam_n \ell x_0).
$$
The sum over $\ell$ is  $\frac 1L (\frac{\sin(\pi L\lam_n x_0)}
{\sin(\pi \lam_n x_0)})^2$ which is always non-negative.  Further if
$n \in {\Cal M}$ then each term in the sum is at least $\cos (2\pi/6)
=\frac 12$ and so the sum here is at least $L/2$.  Thus we see
that
$$
\sum_{\ell = -L}^{L} k\Big( \frac{\ell}{L}\Big) F_1(\ell x_0)
\ge \frac{L}{4}
\sum_{m \in {\Cal M}} f(m) k\Big(\frac{\lam_N-\lam_m}{\lam_N}\Big)
\ge \frac{L}{8} \sum_{m\in {\Cal M}} f(m),
$$
since $\lam_m \in [\lam_N/2,3\lam_N/2]$ for all $m\in {\Cal M}$.
Since $F_1(\ell x_0)=F_1(-\ell x_0)$ we deduce
that for some $1\le \ell_0 \le L$
$$
F_1(\ell_0 x_0)
\ge \frac{1}{8} \sum_{m\in {\Cal M}} f(m) - \frac{F_1(0)}{L-1}
\ge \frac{1}{8} \sum_{m\in {\Cal M}} f(m) -\frac{1}{L-1} \sum\Sb n\\ \lam_n
\le 2\lam_N\endSb f(n).
$$
Using this in (3) we obtain the first assertion of the Lemma.

Suppose now that $\beta \equiv 0\pmod{2\pi}$.  We start with
$$
\int_{-\infty}^{\infty} 2\lam_N K(2\lam_N u) F(x+u)du =
\sum_{n} f(n) \cos(2\pi \lam_n x) k\Big(\frac{\lam_n}{2\lam_N}\Big)
$$
and letting $F_2(x)$ denote the RHS above, we deduce that
$$
F_2(x) \le \max_{u \in [-X/2,X/2]} F(x+u)
+ \frac{2}{\pi^2 X \lam_N} \sum_n f(n).
$$
We then argue as in the preceding paragraph and obtain
the estimate (2).  The case $\beta\equiv \pi\pmod{2\pi}$ follows
since $\cos(t+\pi)=-\cos(t)$.

\enddemo

\head 3. Proof of Theorem 1\endhead

\noindent Let $X$ be large.
Uniformly in
$X\le x\le X^3$ we have (see (12.4.4) of [10])
$$
\Delta(x)= \frac{x^{\frac 14}}{\pi\sqrt {2}}
\sum_{n\le X^3} \frac{d(n)}{n^{\frac 34}}
\cos\Big(4\pi \sqrt{nx}-\frac{\pi}{4}\Big) + O(X^{\epsilon}).
$$
We will apply the result of \S 2 taking
$f(n)=d(n) n^{-\frac 34}$ if $n\le X^3$ and $f(n)=0$ for larger $n$,
$\lam_n = 2\sqrt{n}$, and $\beta = -\frac \pi{4}$.  Then
for $\sqrt{X} \le x \le X^{\frac 32}$ we have $\Delta(x^2)
=\frac{\sqrt{x}}{\pi \sqrt{2}} F(x) + O(X^{\epsilon})$, so that
it suffices to establish an $\Omega$ result for $F$.

Let $L$, $M$ and $N$ be parameters to be
chosen shortly and suppose that $(6L)^{M+1} \le \sqrt{X}$.
Let ${\Cal M}$ be a set of $M$ integers in $[N/4,9N/4]$.
Then (1) of Lemma 3 shows that there exists a point
$x\in [X/2,X^{\frac 32}]$ such that
$$
\align
|F(x)| &\ge \frac 18\sum_{m \in {\Cal M}} \frac{d(m)}{m^{\frac 34}}
-\frac 1{L-1} \sum_{n\le 4N} \frac{d(n)}{n^{\frac 34}}
-\frac{2}{\pi^2 X\sqrt{N}} \sum_{n\le X^3} \frac{d(n)}{n^{\frac 34}}
\\
&\ge \frac{1}{18N^{\frac 34}} \sum_{m\in {\Cal M}} d(m) +
O\Big( \frac{N^{\frac 14}\log N}{L} + \frac{\log X}{X^{\frac 14}}\Big).
\tag{4}\\
\endalign
$$

Choose $L= (\log_2 X)^{10}$ and let $\lambda$ be a
positive real number (we shall see that $\lambda = 2^{\frac 43}$ optimally).
We take ${\Cal M}$ to be the set of integers in $[N/4,9N/4]$
having exactly $[\lam \log_2 N]$ distinct prime factors.
The cardinality of ${\Cal M}$ is
$$
M \asymp \frac{N}{\log N} \frac{(\log_2 N)^{[\lam \log_2 N]-1}}
{([\lam \log_2 N]-1)!}
\asymp \frac{N}{\sqrt{\log_2 N}} (\log N)^{\lam -1 -\lam\log \lam},
$$
upon using Stirling's formula and Theorem 4 of II.6.1 of G. Tenenbaum [9]
for example.  If we take $N= c \log X (\log_2 X)^{1-\lam +\lam\log \lam}
(\log_3 X)^{-\frac 12}$ for a suitably small positive constant $c$ then
the condition $(6L)^{M+1} \le \sqrt X$ is satisfied.  Upon
noting that each $m\in {\Cal M}$ satisfies $d(n) \ge 2^{[\lam \log_2 N]}
\asymp (\log N)^{\lam \log 2}$ we deduce from (4) that
for some $x \in [X/2,X^{\frac 32}]$
$$
\align
|F(x)| &\gg \frac{M (\log N)^{\lam \log 2} }{N^{\frac 34}}
+ O\Big( \frac{N^{\frac 14}}{(\log_2 X)^9} +1 \Big)
\\
&\gg \frac{(\log X)^{\frac 14}}{(\log_3 X)^{\frac 58}}
(\log_2 X)^{\lam\log 2+\frac 34(\lam -1-\lam\log \lam)}
+ O( (\log X)^{\frac 14} (\log_2 X)^{(1-\lam+\lam\log \lam)/4 -9}).
\\
\endalign
$$
The optimal choice of $\lam$ is $\lam=2^{\frac 43}$ which gives the
omega result for $\Delta(x)$ claimed in Theorem 1.

The proof for $P(x)$ is similar.  By modifying the argument in
Titchmarsh [10; \S 12.4] we obtain that uniformly in $X\le x\le X^{3}$
$$
P(x)= -\frac{x^{\frac 14}}{\pi} \sum_{n\le X^3} \frac{r(n)}{n^{\frac 34}}
\cos \Big(2\pi \sqrt{nx}+\frac{\pi}{4}\Big) + O(X^{\epsilon}).
$$
We now apply the result of \S 2 taking $f(n)=r(n)n^{-\frac 34}$
for $n\le X^3$ and $f(n)=0$ for larger $n$, $\lam_n =\sqrt n$ and
$\beta =\frac {\pi}{4}$.  Then for $\sqrt{X} \le x\le X^\frac 32$
we have $P(x^2) = -\frac{\sqrt x}{\pi} F(x) +O(X^{\epsilon})$
so that it suffices to establish an $\Omega$ result for $F$.
Let $L$, $M$ and $N$ be parameters to be chosen and suppose
$(6L)^{M+1} \le \sqrt{X}$.  Let ${\Cal M}$ be a set of $M$ integers in
$[N/4,9N/4]$.  Then (1) of Lemma 3 shows that there
is a point $x\in [X/2, X^{\frac 32}]$ with
$$
|F(x)| \ge \frac{1}{18 N^{\frac 34}}\sum_{m\in {\Cal M}} r(m)
+ O\Big( \frac{N^{\frac 14}}{L} + \frac{1}{X^{\frac 14}}\Big).  \tag{5}
$$
Choose $L=(\log_2 X)^{10}$ and let $\lam$ be a positive real
number (we shall see that the optimal choice of $\lam$ is $2^{\frac 13}$).
We take ${\Cal M}$ to be the set of integers in $[N/4,9N/4]$
having exactly $[\lam \log_2 N]$ distinct prime factors all
of which are $1 \pmod 4$.  Modifying the arguments in II.6 of
Tenenbaum [9] we see that the cardinality of ${\Cal M}$ is
$$
M \asymp \frac{N}{\log N}
\frac{(\frac 12\log_2 N)^{[\lam \log_2 N]-1}}{([\lam \log_2 N]-1)!}
\asymp \frac{N}{\sqrt{\log_2 N}}
(\log N)^{\lam -1 -\lam \log \lam-\lam \log 2}.
$$
If we let $N=c \log X (\log_2 X)^{1-\lam +\lam \log \lam +\lam\log 2}
(\log_3 X)^{-\frac 12}$ for a suitably small positive constant $c$
then the condition $(6L)^{M+1} \le \sqrt{X}$ is met.  Upon noting
that $r(m)\ge 2^{[\lam \log_2 N]} \asymp (\log_2 X)^{\lam \log 2}$
for all $m\in {\Cal M}$ we obtain from (5) that
for some $x\in [X/2,X^{\frac 32}]$
$$
|F(x)| \gg  \frac{(\log X)^{\frac 14}}{(\log_3 X)^{\frac 58}}
(\log_2 X)^{\frac{\lam \log 2}{4} + \frac 34(\lam-1-\lam\log \lam)}
+O((\log X)^{\frac 14} (\log_2 X)^{(1-\lam +\lam \log(2\lam))/4 -9}).
$$
The optimal choice for $\lam$ is $\lam = 2^{\frac 13}$ which
establishes this case of Theorem 1.

\head 4. Proof of Theorem 2 \endhead

\proclaim{Proposition 4} Let $x\ge 2$ and $N\ge 2$ be real numbers.  Then
for a fixed integer $k\ge 2$
$$
\align
\frac{N^{\frac 1k}}{\sqrt{\pi}} \int_{-\infty}^{\infty} &\Delta_k
(x^k e^{u/x}) e^{-u^2 N^{\frac 2k}} du
=  O(x^{\frac k2 -\frac 35} N^{\frac 12+\epsilon})\\
&+\frac{x^{\frac {k-1}{2}}}{\pi \sqrt{k}} \sum_{n=1}^{\infty}
\frac{d_k(n)}{n^{\frac{k+1}{2k}}} \exp(-\pi^2 (n/N)^{\frac 2k})
\cos \Big( 2\pi k n^{\frac 1k}x +\frac{k-3}{4}\pi\Big)
.\\
\endalign
$$
\endproclaim

Assuming Proposition 4 we now prove Theorem 2.  We apply
the result of \S 2 taking $f(n)=d_k(n)n^{-\frac{k+1}{2k}}
\exp(-\pi^2 (n/N)^{\frac 2k})$, $\lam_n = kn^{\frac 1k}$ and
$\beta = \frac{k-3}{4}\pi$.  By Proposition 4 it suffices to
establish $\Omega$ results for the corresponding $F(x)$ where
we suppose that $X/2\le x\le X^2$ say. (The error term in Proposition
4 is negligible for our choice of $N$ which will be $O(X^{\epsilon})$.)

We choose $L= (\log_2 X)^{k^3+20}$ and select ${\Cal M}$ to
be the set of integers in $[2^{-k}N,(3/2)^kN]$ containing
exactly $[\lam \log_2 N]$ distinct prime factors; here $\lam$ is
a positive real number which will be optimally chosen as
$k^{\frac{2k}{k+1}}$.  As in \S 3, we see that
the cardinality of ${\Cal M}$ is
$M \asymp N (\log N)^{\lam -1-\lam\log \lam} (\log_2 N)^{-\frac 12}$.
If we choose
$N =c_k \log X (\log_2 X)^{1+\lam\log \lam -\lam}(\log_3 X)^{-\frac 12}$
for a suitably small positive constant $c_k$ then the condition
$(6L)^{M+1}\le X$ is met.  Since $d_k (m) \ge k^{[\lam \log_2 N]}
\asymp (\log_2 X)^{\lam \log k}$ for each $m\in {\Cal M}$,
Lemma 3 then establishes that for some $x \in [X/2,X^2]$ we have
$$
|F(x)| \gg  (\log X)^{\frac{k-1}{2k}}
(\log_2 X)^{\frac{k+1}{2k}(\lam-1-\lam\log \lam)
+\lam \log k} (\log_3 X)^{-\frac 12 -\frac{k-1}{4k}} +
O\Big(\frac{N^{\frac {k-1}{2k}} (\log N)^{k-1}}{L}\Big).
$$
Choosing optimally $\lam = k^{\frac {2k}{k+1}}$ we obtain
the desired omega result for $F(x)$ and hence Theorem 2.
When $k\equiv 3\pmod 8$ then $\beta\equiv 0\pmod {2\pi}$ and
when $k\equiv 7 \pmod 8$ then $\beta \equiv \pi\pmod {2\pi}$,
and so in these cases Lemma 3 leads to the one sided omega
results claimed in Theorem 2.

It remains lastly to prove Proposition 4.  The proof is
based on a standard procedure using Perron's formula,
shifting contours, invoking the functional equation for
$\zeta(s)$, and then applying the method of stationary phase.
One can also extract Proposition 4 from the work of
Hafner [3] (see (3.2.8)).  For the sake of completeness
we supply a proof.

\demo{Proof of Proposition 4} Write $D_k(x)=\sum_{n\le x} d_k(n)$
and consider
$$
\frac{N^{\frac 1k}}{\sqrt{\pi}}
\int_{-\infty}^{\infty} D_k(x^k e^{u/x}) e^{-u^2 N^{\frac 2k}} du.
$$
By Perron's formula this is, for some $c>1$,
$$
\align
&=\frac{N^{\frac 1k}}{\sqrt{\pi}} \int_{-\infty}^{\infty}
\frac{1}{2\pi i} \int_{c-i\infty}^{c+i\infty} \zeta(s)^k x^{ks}
\exp\Big(\frac{us}{x} - u^2 N^{\frac 2k}\Big) \frac{ds}{s} du
\\
&= \frac{1}{2\pi i} \int_{c-i\infty}^{c+i\infty} \zeta(s)^k x^{ks}
\frac{N^{\frac 1k}}{\sqrt{\pi}} \int_{-\infty}^{\infty} \exp\Big(-
\Big(uN^{\frac 1k}-\frac{s}{2N^{\frac 1k}x}\Big)^2
+\frac{s^2}{4N^{\frac 2k}x^2}\Big) du \frac{ds}{s}\\
&= \frac{1}{2\pi i} \int_{c-i\infty}^{c+i\infty} \zeta(s)^k x^{ks}
\exp\Big( \frac{s^2}{4N^{\frac 2k}x^2}\Big)\frac{ds}{s}.\\
\endalign
$$
We move the line of integration above to the
line $a -i\infty$  to $a +i\infty$ where we take $a=-\frac 1{\log x}$.
The pole at $0$ gives an amount $O(1)$ while the pole at $s=1$ contributes
$$
\mathop{\text{Res}}_{s=1}\ \zeta(s)^k  \frac{x^{ks}}{s}
\exp\Big( \frac{s^2}{4N^{\frac 2k}x^2}\Big)
= \frac{N^{\frac 1k}}{\sqrt{\pi}} \int_{-\infty}^{\infty}
e^{-u^2 N^{\frac 2k}}
\Big(\mathop{\text{Res}}_{s=1}\ \zeta(s)^k  \frac{(x^{k}e^{u/x})^{s}}{s}\Big)
du.
$$
We conclude that
$$
\frac{N^{\frac 1k}}{\sqrt{\pi}}
\int_{-\infty}^{\infty} \Delta_k(x^k e^{u/x}) e^{-u^2 N^{\frac 2k}}   du
= \frac{1}{2\pi i}
\int_{a-i\infty}^{a  +i\infty} \zeta(s)^k x^{ks}
\exp\Big( \frac{s^2}{4N^{\frac 2k}x^2}\Big)\frac{ds}{s} +O(1).
$$
We use the functional equation $\zeta(s)=\chi(s)\zeta(1-s)$
where $\chi(s)= 2^{s-1}\pi^s \sec (\pi s/2)/\Gamma(s)$ and expand
$\zeta(1-s)^k = \sum_{n=1}^{\infty} d_k(n) n^{s-1}$.  Then the
above becomes
$$
\sum_{n=1}^{\infty} \frac{d_k(n)}{n}
\frac{1}{2\pi i} \int_{a-i\infty}^{a+i\infty} \chi(s)^{k} (x^k n)^s
\exp\Big( \frac{s^2}{4N^{\frac 2k}x^2}\Big)\frac{ds}{s} + O(1). \tag{6}
$$

Call the integral in (6) above $I_n$. The integral over the line segment
from $a-i$ to $a+i$ gives an amount $\ll n^a \log x$ and
note that the integrand at $a-it$ is the complex conjugate
of the integrand at $a+it$.  Thus
$$
I_n =  \text{Re }\frac{1}{\pi i} \int_{1}^{\infty} \chi(a+it)^k
(x^{k}n)^{a+it} \exp\Big(\frac{a^2 +2ait -t^2}{4N^{\frac 2k}x^2}\Big)
\frac{dt}{t}
+ O(n^a \log x).
$$
Stirling's formula gives that
$\chi(a+it) = (2\pi/t)^{a+it-\frac 12} e^{i(t+\frac \pi{4})} (1+O(1/t))$.
Hence
$$
I_n=\text{Re }\frac{1}{\pi i}
\int_1^{\infty}  (x^{k}n)^{a +it}
\Big(\frac{2\pi}{t}\Big)^{k(a+it-\frac 12)} e^{ik(t+\frac \pi{4})}
\exp\Big(-\frac {t^2}{4N^{\frac 2k}x^2}\Big) \frac{dt}{t}
+ O(n^a x^{\frac k2 -1+\epsilon}\sqrt{N}). \tag{7}
$$

We use the method of stationary phase (which occurs at
$t=2\pi n^{\frac 1k}x$) to evaluate the above integral.
We split the
cases when $|t-2\pi xn^{\frac 1k}| \le (xn^{\frac 1k})^{\frac 35}$
and when $|t-2\pi xn^{\frac 1k}| >(xn^{\frac 1k})^{\frac 35}$.
In the first case (call $y=t-2\pi x n^{\frac 1k}$ so that
$|y|\le (xn^{\frac 1k})^{\frac 35}$ here)
we get by a Taylor expansion
$$
\align
\text{Re }\frac{(xn^{\frac 1k})^{\frac k2 -1}}{2\pi^2 i}
 \int_{|y|\le (xn^{\frac 1k})^{\frac 35} }
\exp\Big(
-\pi^2 \Big(\frac{n}{N}\Big)^{\frac 2k}
&+i \Big( \frac{k\pi}{4} +2\pi k n^{\frac 1k}x -
\frac{k y^2}{4\pi xn^{\frac 1k}}\Big)\Big)
\\
&\times \Big(1+O\Big(\frac{1}{(xn^{\frac 1k})^{\frac 15}}\Big)\Big)dy.\\
\endalign
$$
Using $\int_{|z|\le T} e^{-iz^2} dz = \sqrt{\pi} e^{-i\pi/4} +O(T^{-1})$
we obtain that the above is
$$
\frac{(xn^{\frac 1k})^{\frac {k-1}{2}}}{\pi \sqrt{k}}
\exp\Big(-\pi^2 \Big(\frac nN\Big)^{\frac 2k}\Big)
\Big( \cos\Big(2\pi kn^{\frac 1k}x+ \frac {k-3}{4}\pi\Big)
+ O\Big(\frac{1}{(xn^{\frac 1k})^{\frac 1{10}}}\Big)\Big) .
$$
To handle the second case we note that for any $1\le y\le
2\pi xn^{\frac 1k} - (xn^{\frac 1k})^{\frac 35}$ we have
(see Lemma 4.2 of [10])
$$
\int_1^y \exp\Big( i\Big(t\log (x^k n) + kt +\frac{k\pi}{4} -
kt\log \frac{t}{2\pi}\Big)\Big) dt
\ll \frac{1}{\log (2\pi x n^{\frac 1k}/y)}.
$$
Using this and integration by parts we see that
the integral in (7) over the range $1\le t\le
2\pi xn^{\frac 1k} - (xn^{\frac 1k})^{\frac 35}$ is
$$
\ll (xn^{\frac 1k})^{\frac k2 -\frac 35} \exp(-(n/N)^{\frac 2k})
+ n^a x^{\frac k2 -1} N^{\frac 12}.
$$
The same bound applies to the integral over the range
$t\ge 2\pi xn^{\frac 1k} +(xn^{\frac 1k})^{\frac 35}$.
Putting these estimates together we find that
$$
\align
I_n &= \frac{(xn^{\frac 1k})^{\frac {k-1}{2}}}{\pi \sqrt{k}}
\exp\Big(-\pi^2 \Big(\frac nN\Big)^{\frac 2k}\Big)
 \cos\Big(2\pi kn^{\frac 1k}x+ \frac {k-3}{4}\pi\Big)
\\
&\hskip 1 in + O(n^{a} x^{\frac k2 -1+\epsilon}\sqrt{N}
+ (xn^{\frac 1k})^{\frac k2-\frac 35} \exp(-(n/N)^\frac{2}{k})).
\\
\endalign
$$
Using this in (6) we obtain the Proposition.

\enddemo

\enddocument